\documentclass [11pt]{amsart} 
\usepackage{epic,eepic,latexsym, amssymb, amscd, amsfonts, color}

\input xy
\xyoption {all}
\def \kv {{\mathsf K_v}}

\def \mv {{\mathsf M_v}}
\def \mw {{\mathsf M_w}}
\def \mvp {{\mathsf M^{+}_{v}}}  
\def \mvm {{\mathsf M^{-}_{v}}}
\def \mwp {{\mathsf M^{+}_{w}}}
\def \mwm {{\mathsf M^{-}_{w}}}
\def \fm {{{\bf R}{\mathcal S}}}
\def \ms {{\mathfrak M}}

\newtheorem{theorem}{Theorem} 
\newtheorem {lemma}{Lemma} 
 
\newtheorem {proposition}{Proposition}

\newcommand{\comment}[1]{}

\theoremstyle{definition}
\newtheorem*{remark}{Remark} 

\begin{document}

\baselineskip=16.75pt
\title[On strange duality for abelian surfaces II]{On the strange duality conjecture for abelian surfaces II}
\author{Barbara Bolognese}
\address{Department of Mathematics, Northeastern University}
\email{bolognese.b@husky.neu.edu}
\author {Alina Marian}
\address {Department of Mathematics, Northeastern University}
\email {a.marian@neu.edu}
\author {Dragos Oprea}
\address {Department of Mathematics, University of California, San Diego}
\email {doprea@math.ucsd.edu}
\author {Kota Yoshioka}
\address{Department of Mathematics, Kobe University}
\email {yoshioka@math.kobe-u.ac.jp}

\begin{abstract} In the prequel to this paper, two versions of Le Potier's strange duality conjecture for sheaves over abelian surfaces were studied. A third version is considered here. In the current setup, the isomorphism involves moduli spaces of sheaves with fixed determinant and fixed determinant of the Fourier-Mukai transform on one side, and moduli spaces where both determinants vary, on the other side. We first establish the isomorphism in rank one using the representation theory of Heisenberg groups. For product abelian surfaces, the isomorphism is then shown to hold for sheaves with fiber degree $1$ via Fourier-Mukai techniques. By degeneration to product geometries, the duality is obtained generically for a large number of numerical types. Finally, it is shown in great generality that the Verlinde sheaves encoding the variation of the spaces of theta functions are locally free over moduli. \end{abstract}
\maketitle

\section{Introduction}

Three versions of Le Potier's strange duality conjecture were formulated in \cite{abelian1} for a polarized abelian surface $(X, H)$. We recall them briefly. 

For a sheaf $E \to X$, we write $$v(E) = \text{ch}\, E \in H^{2\star} (X, \, {\mathbb Z})$$ for its Mukai vector. 
For two Mukai vectors $$v=(v_0, v_2, v_4), \,\,\,w=(w_0, w_2, w_4)\in H^{2\star}(X, \mathbb Z),$$ the Mukai pairing is given by $$\langle v, w\rangle = \int_{X} v_2 w_2-v_0w_4-v_4w_0.$$ We also set standardly $$v^{\vee}=(v_0, -v_2, v_4)\in H^{2\star}(X, \mathbb Z).$$

Let $\mathsf M_v$ be the moduli space of Gieseker $H$-semistable sheaves with Mukai vector $v$. When $v$ is primitive and the polarization $H$ is generic, the moduli space $\mathsf M_v$ consists of stable sheaves only, and is smooth projective of dimension
$$\dim \mv = 2 d_v +2, \, \, \, \text{where} \, \, d_v = \frac{1}{2} \langle v, \, v\rangle .$$ We will make this assumption about the moduli spaces throughout the paper, unless specified otherwise. 
We furthermore consider three subspaces of $\mv$:
\noindent
\begin{itemize}
\item [-] the space 
$\mvp$ of sheaves with a fixed determinant line bundle;
\item [-] the space 
$\mvm$ of sheaves with fixed determinant of their Fourier-Mukai transform;
\item [-] the space 
$\kv$ of sheaves for which both the determinant and the determinant of their Fourier-Mukai transform is fixed.
\end{itemize}
In introducing the spaces $\mvm, \, \kv$, we use the Fourier-Mukai transform $$\fm: {\mathbf D} (X) \longrightarrow {\mathbf D} (\widehat{X})$$ with respect to the standardly normalized Poincar\'{e} line bundle $${\mathcal P} \to X \times \widehat{X}.$$ The moduli space $\mathsf K_v$ is precisely the fiber of the Albanese map $$\mathsf a: \mathsf M_v \to X\times \widehat X.$$ The morphism $\mathsf a$ is defined up to the choice of a reference sheaf $E_0$ of type $v$. Specifically, $$\mathsf a(E)=(\det \fm (E)\otimes \det \fm (E_0)^{\vee}, \det E\otimes \det E_0^{\vee}).$$

Consider now two Mukai vectors $v$ and $w$, orthogonal in the sense that $$\langle v^{\vee}, w\rangle=-\chi(X, v\cdot w)=0.$$ A sheaf $F \to X$ with Mukai vector $$w = \text{ch} (F) \in H^{2\star} (X, {\mathbb Z})$$ gives rise to a line bundle $$ \Theta_F \to \mv$$  by the standard determinant construction described in \cite{lepotier}, \cite{jun}. Specifically, if a universal family $\mathcal E\to \mathsf M_v\times X$ exists, we set \begin{equation}\label{thetamap}\Theta_F=\det {\mathbf R}p_{!}(\mathcal E\otimes q^{\star} F)^{-1}\to \mathsf M_v,\end{equation} where $p, q$ are the two projections. By restriction, one gets line bundles on each of the subspaces $\mvp, \, \mvm, \, \kv.$ 

Within a fixed Mukai class $w$, for each of the four moduli spaces considered, the dependence of the determinant line bundle on $F$ takes a particular form, as explained in \cite{abelian1}:
\begin{itemize}
\item[-] on $\kv$, the line bundle $\Theta_F = \Theta_w$ depends only on the Mukai class $w$ of $F$;
 \item [-] on $\mvp$,  the line bundle $\Theta_F$ is constant as long as the determinant of $F$ is fixed;
 \item[-] on $\mvm,$ the line bundle $\Theta_F$ is constant as long the determinant of the Fourier-Mukai transform of $F$ is fixed;
 \item [-] on $\mv$, the line bundle $\Theta_F$ is constant as long as $F$ has both its determinant and its FM-transform determinant fixed. 
 \end{itemize}

Keeping these variations in mind, we write $\Theta_w$ for the determinant line bundle on each of the four moduli spaces, suitably understood. The distinctions above are further highlighted by the numerical equalities, cf. \cite{abelian1}:
\begin{equation}\label{h1}\chi (\kv, \, \Theta_w) = \chi (\mw, \, \Theta_v) = \frac{d_{v}^{2}}{d_{v}+d_{w}}\binom{d_{v}+d_{w}}{d_{v}},\end{equation}
$$\chi (\mvp, \, \Theta_w) = \chi (\mwp, \, \Theta_v) = \frac{1}{2}\frac{c_{1}(v\otimes
w)^{2}}{d_{v}+d_{w}}\, \binom{d_{v}+d_{w}}{d_{v}},$$
$$\chi (\mvm, \, \Theta_w) = \chi (\mwm, \, \Theta_v) = \frac{1}{2}\frac{c_{1}(\hat v\otimes \hat w)^{2}}{d_{v}+d_{w}}\binom{d_{v}+d_{w}}{d_{v}}.$$ Here, $\widehat v$ and $\widehat w$ denote the cohomological Fourier-Mukai transforms of $v$ and $w$. 
\vskip.1in
\noindent

Following Le Potier's original strange duality proposal \cite{lepotier0}, it was shown in \cite{abelian2} that the Brill-Noether divisors $$\Theta^{+} = \{ (E, F) \, \text{with} \, \, \mathbb H^0 (E \otimes^{\mathbf L} F) \neq 0 \} \subset \mvp \times \mwp$$ and 
$$\Theta^{-} = \{ (E, F) \, \text{with} \, \, \mathbb H^0 (E \otimes^{\mathbf L} F) \neq 0 \} \subset \mvm \times \mwm$$ induce isomorphisms of spaces of sections
$${\mathsf D}^{+}: H^0 (\mvp, \, \Theta_w)^{\vee} \longrightarrow H^0 (\mwp, \, \Theta_v), $$
$${\mathsf D}^{-}: H^0 (\mvm, \, \Theta_w)^{\vee} \longrightarrow H^0 (\mwm, \, \Theta_v), $$
for infinitely many Mukai vectors $v$ and $w$ and for an abelian surface $(X, H)$ which is a product of elliptic curves.

In this paper we focus on the third possible geometry, associated with the divisor 
$$\Theta = \{ (E, F) \, \text{with} \, \, \mathbb H^0 (E \otimes^{\mathbf L} F) \neq 0 \} \subset \kv \times \mw.$$ The current setting is particularly interesting since it exhibits the fixed versus unfixed determinant asymmetry also present for moduli spaces of bundles over curves \cite{beauville}. In this asymmetric setup, 
we establish the duality generically for a large class of Mukai vectors $v$ and $w$, as captured in our main Theorem \ref{t1} below. We now explain the salient points of the argument and state the most important results along the way. 
\vskip.1in

The starting point is the case when $v$ and $w$ are Mukai vectors of rank 1. For each integer $a> 0$, we let $X^{[a]}$ be the Hilbert scheme of $a$ points on $X$, and let $$K^{[a]} \subset X^{[a]}$$ be the generalized Kummer variety of $a$ points adding to zero on $X$. When rank $v$ = rank $w$ = 1, we have $$\kv \simeq K^{[a]}, \,\, \, \mw \simeq X^{[b]} \times \widehat{X}, $$ for suitable $a, b.$ In this setup, we prove

\begin{theorem}\label{t2}
Let $L \to X$ be an ample line bundle on an arbitrary abelian surface. Write $\chi (X, L) = \chi= a+b$ for positive integers $a$ and $b$. The divisor
$$\Theta_L =\{(I_Z, I_W, y) \, \, \text{with} \, \, H^0 (I_Z \otimes I_W \otimes y \otimes L) \neq 0\} \subset K^{[a]} \times X^{[b]} \times \widehat{X}$$
induces an isomorphism
$$D_L: \, H^0 (K^{[a]}, \, \Theta_v)^{\vee} \longrightarrow H^0 (X^{[b]} \times \widehat{X}, \, \Theta_w).$$

\end{theorem}
The analogous isomorphism when both sides involve the Hilbert schemes $X^{[a]}$ and $X^{[b]}$ and the theta bundles over them was shown to hold for all surfaces in \cite {abelian3}. By contrast, Theorem \ref{t2} is a subtler statement specific to abelian surfaces. Its proof requires new ideas and is obtained using the representation theory of the Heisenberg group. 

Paralleling \cite {generic} and \cite{abelian2}, the above result implies  strange duality for product abelian surfaces via Fourier-Mukai techniques. Specifically, for moduli spaces of sheaves which are stable with respect to a suitable polarization in the sense of Friedman \cite {F}, we show

\begin {theorem} \label{t3}Let $X=B\times F$ be a product abelian surface. Assume $v$ and $w$ are two orthogonal Mukai vectors of ranks $r, r'\geq 2$ with $$c_1(v)\cdot f=c_1(w)\cdot f=1.$$ 
Then, the locus $$\Theta = \{ (E, F) \, \, \text{with} \, \, {\mathbb H}^0 (E \otimes^{\mathbf L} F) \neq 0 \} \subset \kv \times \mw$$ is a divisor, and induces an isomorphism  
$$\mathsf D:H^0(\mathsf K_v, \Theta_w)^{\vee}\to H^0(\mathsf M_w, \Theta_v).$$ 
\end {theorem} 

In order to move from the product geometry of Theorem \ref{t3} to a generic abelian surface, we study the Verlinde sheaves $$\mathbb V, \mathbb W\to \mathcal A.$$ These are defined in Section \ref{deg}, and encode the spaces of generalized theta functions  $H^0 (\kv, \, \Theta_w)$ and $H^0 (\mw, \, \Theta_v)$ respectively, as the pair $(X, H)$ varies in its moduli space $\mathcal A$. We need to ensure that the Verlinde sheaves are generically locally free of expected rank given by the holomorphic Euler characteristics \eqref{h1}: $$\text{rank }\mathbb V=\text{rank } \mathbb W=\frac{d_v^2}{d_v+d_w}\binom{d_v+d_w}{d_v}.$$ We establish this in our situation by showing that for surfaces of N\'eron-Severi rank $1$ the theta line bundles are big and nef, and therefore carry no higher cohomology. This yields the following generic strange duality statement, which constitutes our main result. 

\begin{theorem}\label{t1}
Assume $(X, H)$ is a generic primitively polarized abelian surface, and $v, \,  w$ are two orthogonal Mukai vectors of ranks $r, \, r'\geq 2$ with 
\begin {itemize}
\item [(i)] $c_1(v)=c_1(w)=H;$
\item [(ii)] $\chi(v)< 0, \chi(w)< 0$.
\end {itemize}
Then, the locus $$\Theta = \{ (E, F) \, \, \text{with} \, \, {\mathbb H}^0 (E \otimes^{\mathbf L} F) \neq 0 \} \subset \kv \times \mw$$ is a divisor, and induces an isomorphism 
$${\mathsf D}: H^0 (\kv, \, \Theta_w)^{\vee} \longrightarrow H^0 (\mw, \, \Theta_v). $$
\end{theorem}

While the statements of Theorems \ref{t3} and \ref{t1} mirror the $K3$ and abelian cases studied in \cite{generic} and \cite {abelian2}, different arguments are needed in the current {\it asymmetric} abelian setup. Several technical assumptions present in \cite{generic} and \cite{abelian2} are in addition removed, yielding stronger results. 

\vskip.1in
Finally, in Section \ref{locallyfree} we show in great generality that the Verlinde sheaves $$\mathbb V, \mathbb W\to \mathcal A$$ are in fact locally free over the entire moduli space $\mathcal A$ even though the higher cohomology of theta line bundles may not vanish. Specifically, this is implied by the following 

\begin{theorem} \label{t4}
Let $(X, H)$ be a polarized abelian surface. Assume that $$v=(r, dH, \chi), \,\,w=(r', d'H, \chi')$$ are orthogonal primitive Mukai vectors of ranks $r, r'\geq2$ such that \begin {itemize}
\item [(i)] $d, d'>0$;
\item [(ii)] $\chi<0,\,\, \chi'<0.$
\end {itemize} Assume furthermore that if $(d,\chi)=(1, -1)$, then $(X, H)$ is not a product of two elliptic curves. We have 
$$h^0(\mathsf K_v, \Theta_w)=\chi(\mathsf K_v, \Theta_w)=\frac{d_v^2}{d_v+d_w}\binom{d_v+d_w}{d_v}.$$ 
Moreover, for any representative $F \in \mathsf K_w,$
$$h^0(\mathsf M_v, \Theta_F)=\chi(\mathsf M_v, \Theta_F)=\frac{d_w^2}{d_v+d_w}\binom{d_v+d_w}{d_v}.$$ \end{theorem}

The proof uses Bridgeland stability conditions, and relies on recent results concerning wall-crossing as stability varies. As walls are crossed, the dimensions of the space of sections do not change. Crucially, we show that we can move away from the Gieseker chamber to a chamber for which the theta line bundles become big and nef. In order to control the wall-crossings and complete the argument, we make use of the explicit description of the movable cone of the moduli space recently obtained in \cite {Y2}; see also \cite {BM}. 

\subsection{Acknowledgements} 
A.M. and D.O. were partially supported by the NSF grants DMS 1001486, DMS 1150675, DMS 1303389, as well as by Sloan Foundation Fellowships. K.Y. was partially supported by
the Grant-in-aid for Scientific Research 22340010, JSPS. Correspondence with Emanuele Macr\`i is gratefully acknowledged. 

\section{The rank one case}

\subsection{Notation and preliminaries}
We let $X$ be an arbitrary abelian surface and consider two Mukai vectors $v$ and $w$ with $$\text{rank} \, v = \text{rank} \, w = 1.$$ Specifically, letting $L \to X$ be an ample line bundle, and writing $\chi (L) = a+ b$ for positive integers $a$ and $b$,  we set
$$ v = (1,\,  0, \, -a), \,\, \, \, w = (1, \, c_1 (L), a).$$ We then have $$\kv \simeq K^{[a]}, \, \, \, \mw \simeq X^{[b]} \times \widehat X,$$ and the strange duality divisor is 
\begin{equation}
\Theta_L =\{(I_Z, I_W, y) \, \, \text{with} \, \, H^0 (I_Z \otimes I_W \otimes y \otimes L) \neq 0\} \subset K^{[a]} \times X^{[b]} \times \widehat{X}.
\end{equation}
Conforming to standard notation, we next set 
 $$L^{[a]} = \det {\mathbf R}p_{\star} \left ( {\mathcal O}_{\mathcal Z} \otimes q^{\star} L \right ) \, \, \, \text{on} \, \, \, X^{[a]},$$where $ \mathcal Z \subset X^{[a]} \times X$
is the universal subscheme, and $p, q$ are the projections to $X^{[a]}$ and $X$ respectively. Throughout this section we also use  $$L^{[a]} \to K^{[a]}$$ to denote the restriction of the determinant 
line bundle to $K^{[a]} \subset X^{[a]}.$ 

The divisor 
\begin{equation} 
\label{standard}
\Theta^{+}_L = \{(I_Z, I_W) \, \, \, \text{with} \, \, \, H^0 (I_Z \otimes I_W \otimes L) \neq 0 \} \subset X^{[a]} \times X^{[b]},
\end{equation} 
with associated line bundle 
$${\mathcal O} (\Theta^{+}_L ) = L^{[a]} \boxtimes L^{[b]} \, \, \, \text{over}\, \, \, \, X^{[a]} \times X^{[b]}$$
induces an isomorphism
\begin{equation}
\label{classical}
D^{+}_L: H^0 (X^{[a]}, \, L^{[a]} )^{\vee} \longrightarrow H^0 (X^{[b]}, \, L^{[b]} ).
\end{equation}
This constitutes the simplest instance of the strange duality phenomenon on surfaces; the isomorphism is described in \cite{abelian3} and holds uniformly irrespective of the choice of surface. 

Relative to this standard rank one setup, the divisor $\Theta_L$ represents a twist specific to the abelian geometry. In particular, the associated line bundle takes the more complicated form
\begin{equation}
\label{twistedline}
{\mathcal O} (\Theta_L) = L^{[a]} \boxtimes L^{[b]} \boxtimes \widehat{L} \otimes (\mathsf a, \text{id})^{\star} {\mathcal P} \, \, \, \, \, \text{on} \, \, K^{[a]} \times X^{[b]} \times \widehat X,
\end{equation}
where ${\mathcal P} \to X \times \widehat{X}$ is the Poincar\'{e} line bundle, and $$\mathsf a: X^{[b]} \to X$$ denotes the addition of points using the group law. We have also set
$$\widehat L  = \det {\bf R}{\mathcal S} (L)^{-1} \, \, \, \text{on} \, \, \, \widehat{X}. $$ Expression \eqref{twistedline} is obtained by restricting to each factor and using Mumford's see-saw theorem; a detailed explanation is found in Example $1$ of \cite {bundle}. Establishing that the induced map on the spaces of sections 
\begin{equation}
\label{main}
D_L : H^0 (K^{[a]}, \, L^{[a]})^{\vee} \longrightarrow H^0 \left (X^{[b]} \times \widehat{X}, \, L^{[b]} \boxtimes \widehat{L} \otimes ({\mathsf a},\, \text{id}) ^{\star} {\mathcal P} \right )
\end{equation}
is an isomorphism requires new ideas which we now describe.
 
\subsection{Proof of Theorem \ref{t2}} To begin, note that both sides of \eqref{main} have equal dimensions given by the Euler characteristics \eqref{h1}. For the left hand side, this follows from either Lemma $3$ or Example $9$ in \cite{bundle}: both show the vanishing of the higher cohomology of $$L^{[a]}\to K^{[a]}$$ under the assumption that $L\to X$ is ample. For the right hand side, we can invoke Proposition \ref{l2} of Section \ref{locallyfree} which applies to the current context as well. A direct argument is also possible making use of the \'etale pullbacks of the proof below. 

We rephrase the statement of the theorem in two steps. To start, let 
$$\varphi_L: X \longrightarrow \widehat{X}, \, \, \,  \varphi_L (x) = t_{x}^{\star} L \otimes L^{-1}$$ be the Mumford homomorphism; we also make use of 
$ \varphi_{\widehat{L}}: \widehat{X} \longrightarrow X.$ Consider now the diagram
\vskip.2in
\begin{center}
$\xymatrix{K^{[a]} \times X^{[b]} \times \widehat{X} \ar[d]^{\widehat \Phi} \ar[dr]^{\Psi} \\ K^{[a]} \times X^{[b]} \times X \ar[d]^{\Phi} \ar[r]^{\Gamma} & X^{[a]} \times X^{[b]} \\ K^{[a]} \times X^{[b]} \times \widehat{X} 
}$
\end{center}
\vskip.1in
where 
\begin{eqnarray*}
\Phi (I_Z, \, I_W, \, x) &= &\left (I_Z, \, t_{x}^{\star} I_W, \, \varphi_{L} (x)\right ),\\ \widehat \Phi (I_Z, \, I_W, \, y) &=& \left (I_Z, \, I_W, \varphi_{\widehat L} (y) \right ),\\
\Gamma (I_Z, \, I_W, \, x) &= &(t_{-x}^{\star} I_Z, \, I_W), \\ \Psi(I_Z, I_W, y)&=&(t_{-\varphi_{\widehat {L}}(y)}^{\star} I_Z, I_W)\implies \Psi=\Gamma\circ \widehat \Phi.
\end{eqnarray*}
All four maps are \'etale: \begin{itemize} \item [-]$\Phi$ and $\widehat{\Phi}$ have degree $\chi^2 = \chi (L)^2 = \chi (\widehat L)^2$; 
\item [-]$\Gamma$ has degree $a^4$ since it can be viewed as quotienting by the group of $a$-torsion points on $X$;
\item [-]$\Psi = \Gamma \circ \widehat \Phi$ has degree $a^4\chi^2.$
\end{itemize}
We now pull back the divisor $\Theta_L \subset K^{[a]} \times X^{[b]} \times \widehat{X}$ {\it twice}, first by $\Phi$ and then by $\widehat{\Phi}.$ 
\subsubsection{Pullback under $\Phi$}
At the first stage, we obtain
\begin{eqnarray*}
\Phi^{\star} \Theta_L &= &\{ (I_Z, \, I_W, \, x) \, \, \text{with} \, \, H^0 (I_Z  \otimes t_{x}^{\star} I_W \otimes \varphi_L (x) \otimes L) \neq 0\} \\ & = & \{ (I_Z, \, I_W, \, x) \, \, \text{with} \, \, H^0 (I_Z  \otimes t_{x}^{\star} I_W \otimes t_{x}^{\star} L) \neq 0\} \\
& = & \{ (I_Z, \, I_W, \, x) \, \, \text{with} \, \, H^0 (t_{-x}^{\star} I_Z  \otimes  I_W \otimes L) \neq 0\} \\ & = & \Gamma^{\star} \Theta^{+}_L.
\end{eqnarray*}
By contrast with expression \eqref{twistedline}, the line bundle associated with $\Phi^{\star} \Theta_L$ has the simpler form
$${\mathcal O} (\Phi^{\star} \Theta_L) = {\mathcal O}(\Gamma^{\star}\Theta_L^+)=\Gamma^{\star}(L^{[a]}\boxtimes L^{[b]})=L^{[a]} \boxtimes L^{[b]} \boxtimes L^a \, \, \, \text{on} \, \, \, K^{[a]} \times X^{[b]} \times X. $$
The pullback divisor induces the map $\Phi^{\star} D_L$ for which the diagram
\begin{center}
$\xymatrix{D_L: \, \, \, H^0 (K^{[a]}, \, L^{[a]})^{\vee} \ar[r] \ar@{=}[d] & H^0 \left (X^{[b]} \times \widehat{X}, \, L^{[b]} \boxtimes \widehat{L} \otimes ({\mathsf a},\, \text{id}) ^{\star} {\mathcal P} \right ) \ar[d]^{\Phi^{\star}}\\
\Phi^{\star} D_L: H^0 (K^{[a]}, \, L^{[a]})^{\vee} \ar[r] &  H^0 \left (X^{[b]} \times X, \, L^{[b]} \boxtimes {L}^a  \right )
}$
\end{center}
commutes. To show the original duality map $D_L$ is injective (and thus by equality of dimensions an isomorphism), it suffices to show that the simpler $\Phi^{\star} D_L$ in the above diagram is injective.  
\subsubsection{Pullback under $\widehat{\Phi}$}
The second pullback, under $\widehat \Phi$, yields the divisor $$\widetilde \Theta_L = \widehat \Phi^{\star} \Phi^{\star} \Theta_L$$ associated with the line bundle
$${\mathcal O} (\widetilde \Theta_L ) = \widehat \Phi^{\star}(L^{[a]}\boxtimes L^{[b]}\boxtimes L^a)=L^{[a]} \boxtimes L^{[b]} \boxtimes \varphi_{\widehat L}^{\star}L^{a} \, \, \, \text{on} \, \, \, K^{[a]} \times X^{[b]} \times \widehat X. $$ Crucially, by our previous interpretation of $\Phi^{\star} \Theta_L$, we also have
$$\widetilde \Theta_L = \widehat \Phi^{\star} \Phi^{\star} \Theta_L=\widehat \Phi^{\star}\Gamma^{\star}\Theta_L^+=\Psi^{\star} \Theta_L^{+}.$$
By the same argument as before, to show that the original duality map \eqref{main} is an isomorphism, it suffices to show that
\begin{proposition} The morphism
$\widetilde D_L: \, H^0 (K^{[a]}, \, L^{[a]})^{\vee} \longrightarrow  H^0  (X^{[b]} \times \widehat X, \, L^{[b]} \boxtimes \varphi_{\widehat L}^{\star}{L}^{a}) $ induced by $\widetilde \Theta_L$ is injective. 
\end{proposition}
\proof We interpret the duality map representation-theoretically, using the theory of discrete Heisenberg groups. $\widetilde D_L$ is better suited for such an interpretation than the seemingly simpler morphism $\Phi^{\star} D_L$ obtained at the previous stage.

We have seen above that $$\Psi^{\star} \mathcal O(\Theta_L^+)=\Psi^{\star} (L^{[a]} \boxtimes L^{[b]}) = L^{[a]} \boxtimes L^{[b]} \boxtimes \varphi_{\widehat L}^{\star} L^{a}.$$ Up to numerical equivalence on $\widehat X$, we have $\varphi_{\widehat L}^{\star}  L=\widehat L^{\chi}$. Thus, there exists $y\in \widehat X$ such that $$\varphi_{\widehat L}^{\star}L= t_{y}^{\star}\widehat L^{\chi}.$$ We define $$M=L\otimes y$$ and calculate $$\widehat M:=\det \fm (M)^{-1}\implies \widehat M=t_y^{\star} \widehat L\implies \widehat M^{a\chi}=t_{y}^{\star}\widehat L^{a\chi}=\varphi_{\widehat L}^{\star} L^{a}.$$ Therefore, $$\Psi^{\star}(L^{[a]}\boxtimes L^{[b]})=L^{[a]}\boxtimes L^{[b]}\boxtimes \widehat M^{a\chi}.$$

We let $\mathsf G({\widehat M}^{\,a} )$ be the Heisenberg group of the line bundle ${\widehat M}^{\,a} \to {\widehat X},$ sitting in an exact sequence
 $$1\rightarrow \mathbb C^{\star} \rightarrow \mathsf G({\widehat M}^{\,a} ) \rightarrow \mathsf H ({\widehat M}^{\,a} )\rightarrow 1,$$
where the quotient is the abelian group $$\mathsf H ({\widehat M}^{a} )= \{y, \, \, t_{y}^{\star} {\widehat M}^{a} \simeq {\widehat M}^{a} \} \subset \widehat X.$$ For an introduction to Heisenberg group actions in the theory of abelian varieties we refer the reader to \cite {Mu}, for instance. 

Importantly, by construction, the \'etale morphism 
$$\Psi: K^{[a]} \times X^{[b]} \times \widehat X \longrightarrow X^{[a]} \times X^{[b]}$$ can be viewed precisely as quotienting by the abelian group $\mathsf H ({\widehat M}^{a} ).$ The latter acts on $K^{[a]} \times \widehat X$ via $$\eta\cdot (I_Z, y)=(t^{\star}_{\varphi_{\widehat M}(\eta)} I_Z, y+\eta)$$ and trivially on $X^{[b]}$. Thus, as a pullback of $L^{[a]} \to X^{[a]}$ under the quotienting map $\Psi$, the line bundle $$L^{[a]} \boxtimes {\widehat M}^{a\chi}\to K^{[a]} \times \widehat X$$ is $\mathsf H({\widehat M}^a)$-equivariant, in other words it is $\mathsf G ({\widehat M}^a )$-equivariant, such that the center acts with weight 0. 
Independently, it is clear that the line bundle  ${\widehat {M}^{a\chi}} \to {\widehat X}$ is $\mathsf G ({\widehat M}^{\,a})$-equivariant, the center acting with weight $\chi$. 
It follows that $$L^{[a]} \to K^{[a]}$$ 
is also $\mathsf G ({\widehat M}^a )$-equivariant, so that the center acts with weight $-\chi$. The spaces of sections $$H^0 (\widehat X, \, {\widehat M}^{a\chi}) \text{ and }H^0 (K^{[a]}, \, L^{[a]})$$ are in turn acted on with weights $\chi$ and $-\chi$ respectively, and furthermore, we can write $$H^0 (X^{[a]}, \, L^{[a]}) = \left ( H^0 (\widehat X, \, {\widehat M}^{a\chi}) \otimes H^0 (K^{[a]}, \, L^{[a]}) \right )^{\mathsf H ({\widehat M}^a )}.$$
Taking into account the long-known isomorphism
$$D_L^{+}: \, H^0 (X^{[a]}, \, L^{[a]} )^{\vee} \longrightarrow H^0 (X^{[b]}, \, L^{[b]})$$ of equation \eqref{classical}, we see that the dual of the linear map $\widetilde D_L$ is
the natural 
$$\widetilde D_L^{\vee}: \, \left ( H^0 (\widehat X, \, {\widehat M}^{a\chi}) \otimes H^0 (K^{[a]}, \, L^{[a]}) \right )^{\mathsf H ({\widehat M}^a )} \otimes H^0 (\widehat X, \, {\widehat M}^{a\chi} )^{\vee} \longrightarrow H^0 (K^{[a]}, \, L^{[a]}),$$
which pairs the vector space $H^0 (\widehat X, \, {\widehat M}^{a\chi} )$ and its dual. To conclude the proposition, we show now that this map is surjective. 

Let $\{\mathsf S_\alpha\}_{\alpha \in I}$ denote the irreducible representations of $\mathsf G({\widehat M}^{a})$ with the center acting with weight $-\chi$. Decomposing into irreducibles, we write
$$H^0 (K^{[a]}, \, L^{[a]}) = \bigoplus_{\alpha} {\mathsf S_{\alpha}} \otimes {\mathbb C}^{m_{\alpha}},\, \, \, 
H^0 (\widehat X, \, {\widehat M}^{a\chi} )^{\vee} =  \bigoplus_{\alpha} {\mathsf S_{\alpha} }\otimes {\mathbb C}^{n_{\alpha}},$$
and the duality map $\widetilde D_L^{\vee}$ is 
$$\widetilde D_L^{\vee}:   \left ( \bigoplus_{\alpha} \left({{\mathbb C}^{n_{\alpha}}}\right)^{\vee}\otimes {\mathbb C}^{m_{\alpha}}\right ) \bigotimes \left ( \bigoplus_{\beta} {\mathsf S_{\beta} }\otimes {\mathbb C}^{n_{\beta}}\right) \longrightarrow 
\bigoplus_{\alpha} {\mathsf S_{\alpha}} \otimes {\mathbb C}^{m_{\alpha}},$$
given explicitly by the natural pairing of the multiplicity spaces $ \left({{\mathbb C}^{n_{\alpha}}}\right)^{\vee}$ and $ {{\mathbb C}^{n_{\alpha}}}.$

We conclude $\widetilde D_L^{\vee}$ fails to be surjective only if there is an irreducible $\mathsf S_{\alpha}$ which appears with nonzero multiplicity $m_{\alpha} \neq 0$ in $H^0 (K^{[a]}, \, L^{[a]})$, but fails to appear in $ H^0 (\widehat X, \, {\widehat M}^{a\chi} ),$
so $n_\alpha = 0.$ This is precluded by the following result, which in level $2$ is Proposition $3.7$ in \cite{I}. This ends the proof of the proposition, and therefore of Theorem \ref{t2}. \qed

\begin{lemma} 
Let $A$ be an abelian surface and $M \to A$ an ample line bundle. For any integer $k \geq 0,$ all irreducible representations with central weight $k$ of the Heisenberg group $\mathsf G (M)$ appear in the $\mathsf G (M)$-module $H^0 (A, \, M^k)$ with nonzero multiplicity.
\end{lemma}

For the benefit of the reader, we give the quick argument, which we lifted from \cite {I}.
Consider the natural homomorphism $\mathsf G(M)\to \mathsf G(M^k)$ and write $$\mathsf K\cong \mathsf G(M)/\mu_k$$ for its image. Fix $\mathsf S$ a representation of the Heisenberg group $\mathsf G(M)$ of weight $k$. Certainly, $\mathsf S$ is a representation of $\mathsf K$ with weight $1$. The induced representation $$\mathsf R=\mathsf {Ind}_{\mathsf K}^{\mathsf G(M^k)}\mathsf S$$ of the Heisenberg group $\mathsf G(M^k)$ has weight $1$, hence it splits as a sum of copies of the unique irreducible representation $H^0(A, M^k)$ of weight $1$: $$\mathsf R=H^0(A, M^k)\oplus\ldots\oplus H^0(A, M^k).$$ We restrict this decomposition to $\mathsf G(M)$. By definition, the induced representation $\mathsf R$ must contain a copy of $\mathsf S$ as a $\mathsf K$-submodule, and therefore also as a $\mathsf G(M)$-submodule. We conclude that $\mathsf S$ must appear in the $\mathsf G(M)$-module $H^0(A, M^k)$, as claimed. 

\section{Product abelian surfaces}

Relying on the rank one case just established, Theorem \ref{t3} is derived by techniques developed in \cite{generic} and \cite{abelian2}.
Specifically, we let $$X = B \times F \to B$$ be a product of elliptic curves, which we view as an abelian surface elliptically fibered over $B$. We write $f$ for the class of the fiber over the origin, and $\sigma$ for the zero section of the fibration. As in \cite {abelian2}, stability of sheaves over $X$ is with respect to a polarization $$H=\sigma+Nf$$ for $N$ large enough. This polarization is suitable in the sense of Friedman \cite {F}.  Assuming $v$ and $w$ are vectors with $$c_1(v)\cdot f=c_1(w)\cdot f=1,$$ 
we show that $$\mathsf D:H^0(\mathsf K_v, \Theta_w)^{\vee}\to H^0(\mathsf M_w, \Theta_v)$$ is an isomorphism. 

As in \cite {abelian2}, we use a fiberwise Fourier-Mukai transform
$$\fm^{\dagger}: \mathbf D(X)\to \mathbf D(X)$$
to move from the rank 1 situation to higher rank Mukai vectors. The kernel of $\fm^{\dagger}$ is given by the pullback of the normalized Poincar\'e sheaf $$\mathcal P_F\to F\times F$$ to the product $X\times_{B} X \cong F \times F\times B.$ The Fourier-Mukai transform gives rise to two birational isomorphisms $$\mathsf K_v\dasharrow K^{[d_v]}$$ and $$\mathsf M_w\dasharrow X^{[d_w]}\times \widehat X$$ which are regular in codimension $1$. Explicitly, for any $E\in \mathsf K_v$ and $F\in \mathsf M_w$, away from codimension two loci, 
Proposition 1 of \cite {abelian2} in conjunction with Theorem 1.1 of \cite {Br1} shows that 
\begin{equation}\label{id1}\fm^{\dagger}(E^{\vee})=I_{Z}(r\sigma-\chi f)[-1],\end{equation}\begin{equation}\label{id2}\fm^{\dagger}(F) = I_{W}^{\vee}  \otimes {\mathcal O} (-r' \sigma +\chi' f) \otimes y^{-1},\end{equation}
for subschemes $$Z\in K^{[d_v]},\,\, W \in X^{[d_w]}, \,\text{ and a line bundle }y\in \widehat X.$$ Here, we wrote 
$$r=\text{rank }(v),\,\, \chi = \chi(v),\,\, r'=\text{rank }(w),\,\,\, \chi'=\chi(w).$$ 
We set
$$L = {\mathcal O} \left (( r+r') \sigma -(\chi+\chi') f \right)\implies \chi(L)=d_v+d_w.$$

Now, the key to finishing the proof is the calculation:
\begin{eqnarray*}
\mathbb H^0 (E \otimes^{\mathbf L} F) &=& \text{Hom}_{{\mathbf D}(X)} (E^{\vee},\, F) =  \text{Hom}_{{\mathbf D}(X)} \left (\fm^{\dagger}(E^{\vee}), \,\fm^{\dagger} (F) \right ) \\
& = & \text{Ext}^1 (I_Z \otimes y \otimes L, \, I_W^{\vee} )= \text{Ext}^1 (I_W^{\vee}, \, I_Z \otimes y \otimes L)^{\vee} \\
& = & \mathbb H^1 (I_W \otimes^{\mathbf L} I_Z \otimes y \otimes  L)^{\vee}.
\end{eqnarray*}
On the locus (of codimension 2 complement) of non-overlapping $(Z, W)$, the last hypercohomology group coincides with the regular cohomology group,
$$\mathbb H^1 (I_W \otimes^{\mathbf L} I_Z \otimes y \otimes  L) =  H^1 (I_W \otimes I_Z \otimes y \otimes L). $$
Thus under the birational map
$${\mathsf K}_v \times {\mathsf M}_w \dasharrow K^{[d_v]} \times X^{[d_w]}\times \widehat X,$$ the two theta divisors $$\Theta = \{(E, \, F): \, \, \mathbb H^0 (E \otimes^{\mathbf L}  F) \neq 0 \} \subset {\mathsf K}_v \times {\mathsf M}_w,$$ and 
$$\Theta_L = \{(I_Z, \, I_W, y): \, \, H^0 (I_Z \otimes I_W \otimes y\otimes L) \neq 0 \}  \subset K^{[d_v]} \times X^{[d_w]}\times \widehat X$$ coincide, and the theta line bundles on each factor match up as well. Since in rank $1$, $\Theta_L$ induces a strange duality isomorphism by Theorem \ref{t2}, the same must be true about the divisor $\Theta$ inducing the map
$$\mathsf D: H^0 ( {\mathsf K}_v, \, \Theta_w)^{\vee} \longrightarrow H^0 ({\mathsf M}_w, \, \Theta_v).$$
This completes the proof. \qed

\begin {remark} The assumption that the rank is at least $3$ is made in \cite{abelian2} to justify that equations \eqref{id1} and \eqref{id2} hold in codimension $1$. This assumption is however not needed, as we now show. The reader wishing to go on to the proof of generic strange duality contained in the next section may choose to skip this argument. 

To begin, we note that identity \eqref{id2} follows from \eqref{id1} via Grothendieck duality. In turn, equation \eqref{id1} is a consequence of the fact that $\fm^{\dagger}(E^{\vee})[1]$ is torsion free, cf. Proposition $1$ in \cite {abelian2}. We will explain that this assertion holds in codimension $1$, in rank $2$. To this end, regard the kernel of $\fm^{\dagger}$, namely the Poincar\'e sheaf $$\mathcal P\to X\times_{B}X,$$ as an object over $X\times X$ via the diagonal embedding $$X\times_{B}X\to X\times X.$$ We will prove
\begin{lemma}
\label{finite}
For all sheaves $E$ away from a codimension 2 locus in the moduli space,  the set  $$T_E=\{x \in X: \text{\upshape{Hom}}(E,{\mathcal P}_{|X \times \{x \}}) \neq 0 \} \subset X$$ is finite.
\end{lemma}
Assuming the lemma, we show that for all $E$ such that $T_E$ is a  finite set, the transform $\fm^{\dagger}(E^{\vee})[1]$ is a torsion free sheaf. To see this, consider a locally free resolution
\begin{equation}\label{res2}0 \to V \to W \to E \to 0\end{equation}
such that $W={\mathcal O}_X(-m H)^{\oplus k}$ for sufficiently large $m$.
Then $$\text{Ext}^1(W,{\mathcal P}_{|X \times \{ x\}})=\text{Ext}^2(W,{\mathcal P}_{|X \times \{ x\}})=0,$$ for all $x \in X$. As a consequence, the sheaf $$\widehat W:=\fm^{\dagger}(W^{\vee})$$ is locally free. 
Next, $$\text{Ext}^2(E,{\mathcal P}_{|X \times \{ x\}})=
\text{Hom}({\mathcal P}_{|X \times \{ x\}}, E)^{\vee}=0,$$ using that $E$ is torsion free and $\mathcal P|_{X\times \{x\}}$ is of rank $0$. From the exact sequence induced by the resolution \eqref{res2}, we conclude that $$\text{Ext}^1(V,{\mathcal P}_{|X \times \{ x\}})=\text{Ext}^2(V,{\mathcal P}_{|X \times \{ x\}})=0$$ for all $x\in X$. Therefore, $$\widehat{V}:=\fm^{\dagger}(V^{\vee})$$ is locally free as well.  
The same resolution also shows that we have an exact triangle $$
\fm^{\dagger}(E^{\vee}) \to \widehat{W} \to \widehat{V} \to\fm^{\dagger}(E^{\vee})[1]$$ which induces an exact sequence in cohomology sheaves
$$
0 \to \mathcal H^0(\fm^{\dagger}(E^{\vee})) \to
\widehat{W} \overset{\phi}{\to} \widehat{V} \to 
\mathcal H^1(\fm^{\dagger}(E^{\vee})) \to 0.
$$ 
Note that $\phi_{|\{x \}}$ is injective whenever $x \not \in T_E$. 
Then our assumption implies that $\phi$ is injective
as a morphism of sheaves. Furthermore, $\text{Coker } \phi$ is  torsion free, as claimed. \qed

\vskip.1in

\noindent
{\it Proof of Lemma \ref{finite}.} Consider the set $$\Sigma=\{E: \text{ there exists a fiber } f \text{ such that } E|_{f} \text { contains a subbundle of slope }>1\}.$$ This set has codimension at least $2$ in the moduli space by Lemma $5.4$ of \cite {BH}. (A shift by $1$ in the slope is necessary to align with the numerical conventions of \cite {BH}.) We will assume that $E$ is chosen outside $\Sigma$. Furthermore, we may assume that there is at most one point of the surface where $E$ fails to be locally free. This is always true in the moduli space away from codimension $2$. 

We claim that in this situation $T_E$ consists of finitely many points. Indeed, let $x\in T_E$. Three cases need to be considered. 
\begin {itemize}
\item [(a)] First, we rely on the fact that the polarization is suitable. In this case, the restriction of $E$ to a generic fiber is stable. If $x$ lies on such a generic fiber, then as a consequence of stability, we obtain the vanishing $$\text{Hom}(E,{\mathcal P}_{|X \times \{x \}}) =0.$$ Therefore in this case $x\not \in T_E$.

\item [(b)] Assume now that $x$ lies on a fiber $f$ over which the restriction of $E$ is locally free but unstable. In this situation, $E|_{f}$ splits as $$E|_{f}=S_0\oplus S_1$$ where $S_0$ is a degree zero line bundle over $f$, while $S_1$ has degree $1$. Any other splitting type is not allowed by the definition of $\Sigma$. Now, $$\text{Hom}(E, \mathcal P|_{X\times \{x\}})=\text{Hom}(S_0, \mathcal P|_{X\times \{x\}})\neq 0\implies S_0=\mathcal P|_{X\times \{x\}}.$$ This shows that $x$ must be the point corresponding to the line bundle $S_0$. Since by (a), there are only finitely many unstable fibers, we conclude that there are only finitely many choices for $x$. 

\item [(c)] Finally, we analyze the case when $x$ lies on a fiber over which $E$ is not locally free. Let $s$ be the unique point where $E$ fails to be locally free, and let $f_s$ be the fiber through $s$. Then $$E|_{f_s}=\mathbb C_{s}\oplus F,$$ where $F$ is a rank $2$ degree $0$ vector bundle over $f_s$. If $F$ is semistable, there exists an extension $$0\to S \to F\to S\to 0$$ where $S$ is a line bundle of degree $0$ over $f_s$. We have $$\text{Hom}(E, \mathcal P|_{X\times \{x\}})=\text{Hom}(F, \mathcal P|_{X\times \{x\}})\neq 0\implies \text{Hom}(S, \mathcal P|_{X\times \{x\}})\neq 0$$ $$\implies S=\mathcal P|_{X\times \{x\}}.$$ This proves that $x$ is the point of the fiber through $s$ corresponding to $S$. 

To complete the argument, it suffices to show that the situation when $F$ is not semistable corresponds to a codimension $2$ subset of the moduli space. To this end, consider the codimension $1$ locus $\mathcal Z$ of sheaves in the moduli space which fail to be locally free at exactly one point. This is an irreducible subset. Indeed, any sheaf in $\mathcal Z$ sits in an exact sequence $$0\to E\to E^{\vee\vee}\to \mathbb C_{s}\to 0,$$ with $M=E^{\vee\vee}$ stable locally free of Mukai vector $$v^{\vee\vee}=v+(0, 0, 1).$$ Letting $\mathcal M$ denote the moduli space of such locally free sheaves, there exists a fibration $$\pi:\mathcal Z\to \mathcal M$$ whose fibers over $M$ are Quot schemes of length $1$ quotients $q:M\to \mathbb C_s\to 0$. The sheaf $E$ is recovered uniquely as the kernel of the pair $(M, q)$. Since the fibers of $\pi$ are irreducible of dimension $3$, $\mathcal Z$ must be irreducible as well. 

Now, for locally free sheaves $M\in \mathcal M$, there are finitely many fibers for which $M|_{f}$ is unstable. Consider $$\mathcal Z^{\circ}\hookrightarrow \mathcal Z$$ the set of pairs $(M, q: M\to \mathbb C_s\to 0)$ where $s$ does not lie on an unstable fiber.  The restriction of $M|_{f_s}$ is the Atiyah bundle of rank $2$ and degree $1$. The kernel of $q$ is a torsion free sheaf $E$ which is not locally free at $s$. In fact, we calculate $$E|_{f_s}=\mathbb C_s\oplus F$$ where $F$ is a subsheaf of degree $0$ of the Atiyah bundle $M|_{f_s}$. Since $M|_{f_s}$ is stable, all its proper subbundles have slope $\leq 0$. It follows that $F$ is semistable. Thus, to get $F$'s which are not semistable, we need to select $(M, q)$ from $\mathcal Z\setminus \mathcal Z^{\circ}$. Clearly, $$\mathcal Z\setminus \mathcal Z^{\circ}\to \mathcal M$$ has projective fibers of dimension $2$. Thus, $\mathcal Z\setminus \mathcal Z^\circ$ has codimension $1$ in $\mathcal Z$, as claimed. This completes the proof of Lemma \ref{finite} and ends the remark. 

\end {itemize} 
\qed
\end {remark} 

\section{Generic strange duality}\label{deg}
The isomorphism we established for product abelian surfaces implies strange duality for generic abelian surfaces. This is achieved via degeneration; see also Section $3$ of \cite {generic}. 

Specifically, we let $\mathcal A$ denote the moduli stack of pairs $(X, H)$ with $H^2=2n$, where $H$ is a primitive ample line bundle over $X$. Consider the universal family $$\pi:(\mathcal X, \mathcal H)\to \mathcal A.$$ Fix integers $\chi, \chi'$ and ranks $r, r'\geq 2$. For each $t\in \mathcal A$ representing a polarized abelian surface $(\mathcal X_t, \mathcal H_t)$, consider two orthogonal Mukai vectors $$v_t=(r, c_1(\mathcal H_t), \chi), \,\,\, w_t=(r', c_1(\mathcal H_t), \chi').$$ We form the relative moduli spaces of $\mathcal H_t$-semistable sheaves of type $v_t$ and $w_t$$$\pi:\mathfrak K[v]\to \mathcal A, \,\,\,\,\pi:\mathfrak M[w]\to \mathcal A.$$ The product $$\pi:\mathfrak K[v]\times_{\mathcal A} \mathfrak M[w]\to \mathcal A$$ carries the relative Brill-Noether locus $$\Theta[v,w]=\{(X, H, E, F): \, \mathbb H^0(X, E\otimes^{\mathbf L} F)\neq 0\}$$ obtained as the vanishing of a section of the relative theta line bundle $$\Theta[w]\boxtimes \Theta[v]\to \mathfrak K[v]\times_{\mathcal A} \mathfrak M[w].$$ Pushing forward to $\mathcal A$ via the natural projections $\pi$, we obtain the sheaves $$\mathbb V=\pi_{\star}\left(\Theta[w]\right),\,\, \mathbb W=\pi_{\star}\left(\Theta[v]\right),$$ as well as a section $\mathsf D$ of $\mathbb V\otimes \mathbb W.$ The constructions are explained in detail in \cite {fourier}.

Crucial to the specialization procedure which yields generic strange duality is the statement that $\mathbb V$ and $\mathbb W$ are generically vector bundles of equal rank $$\frac{d_v^2}{d_v+d_w}\binom{d_v+d_w}{d_v}$$ whose fibers are the spaces of generalized theta functions. This is established in Proposition \ref{bignef} below. Assuming this result, we let $\mathcal A^{\circ}\hookrightarrow \mathcal A$ denote the maximal open locus where the generic rank is achieved. Consider also the Humbert locus $$\mathcal S\hookrightarrow \mathcal A$$ of split abelian surfaces $$(X, H)=(B\times F, L_B\boxtimes L_F),$$ for line bundles $L_B\to B, L_F\to F$ of degrees $1$ and $n$. Just as in Section $3$ of \cite {generic}, Theorem \ref{t3} can be rephrased as the statement that $$\mathcal S\hookrightarrow \mathcal A^{\circ}$$ and that furthermore $$\mathsf D:\mathbb V^{\vee}\to \mathbb W$$ is an isomorphism along $\mathcal S$. To make the above claim, we need to exchange stability relative to a suitable polarization required by Theorem \ref{t3} with stability relative to the polarization $H$ (which may lie on a wall). The next section, in particular Proposition \ref{wallc}, shows that the ensuing moduli spaces agree in codimension $1$. We need to pass to the moduli stacks to invoke the proposition, but the corresponding spaces of sections do not change, as explained in Section 3 of \cite {generic}. 

As a consequence, $\mathsf D$ is an isomorphism generically over $\mathcal A^\circ$. Since the generic fibers of $\mathbb V$ and $\mathbb W$ over $\mathcal A^{\circ}$ are spaces of generalized theta functions, we conclude that generic strange duality holds as in Theorem \ref{t1}.\qed \vskip.1in

We now turn to Proposition \ref{bignef} which was used in the argument above. A general local-freeness statement for the Verlinde sheaves will be proven in Section \ref{locallyfree}, but in its context, the proposition gives stronger positivity results with a simpler proof. We show
  
\begin {proposition} \label{bignef} Let $X$ be an abelian surface of Picard rank $1$, with $H$ the generator of the N\'{e}ron-Severi group of $X$. Let $$v=(r, H, \chi), \,\,\,w=(r', d'H, \chi')$$ be two orthogonal vectors of positive rank such that $\chi\neq 0$, $\chi'\leq 0.$ Then, for any $F\in \mathsf K_w$, the line bundle $$\Theta_w:=\Theta_F\to \mathsf M_v$$ is big and nef, hence without higher cohomology. If $\chi'<0$, then the above line bundle is ample. By restriction, the same results hold for $\Theta_w\to \mathsf K_v$. 
\end {proposition}
 
\proof In the $K3$ case, reflections along rigid sheaves were used to conclude that $\Theta_w\to \mathsf M_v$  is big and nef, hence without higher cohomology, cf. Proposition $4$ of \cite {generic}. Unlike $K3$ surfaces, abelian surfaces do not admit rigid sheaves. A different argument will be given.  

The starting point is the following well-known result of Jun Li \cite {junli2}. Specifically, setting $$w_0=(0, rH, -2n),$$ the line bundle $\Theta_{w_0}\to \mathsf M_v$ is big and nef. We will moreover show that for the vector $$w_1=(2n, -\chi H, 0),$$ the line bundle $\Theta_{w_1}\to \mathsf M_v$ is also big and nef. Since for $\chi(w)\leq 0$, $w$ is a linear combination with non-negative coefficients of $w_0$ and $w_1$, the conclusion follows. 

To prove the claim about $w_1$, we consider two cases depending on the sign of $\chi(v)$. Let us first assume that $\chi(v)<0$. By Proposition $3.5$ of \cite {Y}, the shifted Fourier-Mukai transform $\Phi$ with kernel $$\mathcal P[1]\to X\times \widehat X$$ induces an isomorphism of moduli spaces $$\Phi:\mathsf M_v\simeq \mathsf M_{\widehat v}\,\,\, \text {  where }\widehat v=(-\chi, \widehat H, -r) \text{ is a vector on }\widehat X.$$ For $\widehat w=(0, -\chi \widehat H, -2n)$, the bundle $$\Theta_{\widehat w}\to \mathsf M_{\widehat v}$$ is big and nef, again by Jun Li's result. To conclude, it remains to observe that $$\Phi^{\star} \Theta_{\widehat w}=\Theta_{w_1},$$ hence the latter line bundle is also big and nef. 

When $\chi(v)> 0$, the argument is similar. By Proposition $3.2$ of \cite {Y}, we have an isomorphism $$\Psi:\mathsf M_v\simeq \mathsf M_{\widehat v},\,\,\,\, \widehat v=(\chi, \widehat H, r)$$ induced by the composition of the Fourier-Mukai transform with kernel $\mathcal P$ with the dualization. Under this isomorphism, Jun Li's bundle $\Theta_{\widehat w}$, where $\widehat w=(0, \chi \widehat H, -2n)$, corresponds to $\Theta_{w_1}$. \qed

\section {Variation of polarization for the moduli space of Gieseker sheaves} Let $X$ be an arbitrary abelian surface, and fix a Mukai vector $$v:=(r,\xi, a) \in H^*(X,{\mathbb Z})$$ with
$r>0$. For an ample divisor $H$ on $X$, denote by  $$\ms(v), \,\,\ms_H(v)^{ss} \text{ and } \ms_H(v)^{\mu-ss}$$ the stacks of all sheaves, of Gieseker $H$-semistable sheaves, and of slope $H$-semistable  sheaves respectively -- all of type $v$. 

We are concerned with moduli spaces of sheaves when Gieseker stability varies: we show that they agree in codimension $1$ each time a wall is crossed. This fact was used in the degeneration argument of Section \ref{deg} to exchange the suitable polarization with the polarization determined by the first Chern class. 

First, for generic polarizations, the dimension of the moduli space is given by the following Lemma 4.3.2 in \cite {MMY2}:
\begin{lemma}
If $H$ is general with respect to $v$, that is,
$H$ does not lie on a wall with respect to $v$,
then
\begin{equation}
\dim \ms_H(v)^{ss}=
\begin{cases}
\langle v,v \rangle+1,& \langle v, v \rangle>0\\
\langle v,v \rangle+\ell,& \langle v, v \rangle=0\\
\end{cases},
\end{equation} 
where
$\ell=\gcd(r,\xi,a)$.
\end{lemma}

For the purposes of Section \ref{deg}, we also need to analyze the situation when the polarization may lie on a wall. To this end, let $H_1$ be an ample divisor on $X$ which belongs to
a wall $W$ with respect to $v$ and $H$ an ample
divisor which belongs to an adjacent chamber. 
Then Gieseker $H$-semistable sheaves are slope $H_1$-semistable $$\ms_{H}(v)^{ss}\hookrightarrow \ms_{H_1}^{ss}(v)\hookrightarrow \ms_{H_1}(v)^{\mu\text{-}ss}.$$ All these stacks have dimension $\langle v, v\rangle +1$ by Lemma $3.8$ of \cite{KY}. 
We estimate 
the codimension of 
$$\ms_{H_1}(v)^{ss} \setminus \ms_H(v)^{ss}.$$
Specifically, we prove
\begin{proposition}\label{wallc} Assume that $v$ is a Mukai vector of positive rank with the property that there are no isotropic vectors $u$ of positive rank such that $\langle v, u\rangle =1 \text{ or } 2$. Then, 
\begin{equation}
\begin{split}
(\langle v,v \rangle+1)-
\dim(\ms_{H_1}(v)^{ss} \setminus \ms_H(v)^{ss}) 
\geq 2.
\end{split}
\end{equation}
Therefore, in this situation, 
$\ms_H(v)^{ss}$ 
is independent of the choice of ample line bundle $H$ (generic or on a wall)
away from codimension $2$. 

The same statement holds true for the moduli stack $\mathfrak K_H(v)^{ss}$ of sheaves with fixed determinant and fixed determinant of the Fourier-Mukai. 

\end{proposition}

\proof The proof is essentially contained in Proposition $4.3.4$ of \cite {MMY2}, but since specific aspects of the argument are used below, we give an outline for the benefit of the reader. Let  $E$ be a Gieseker $H_1$-semistable sheaf, which is however not  Gieseker $H$-semistable. In particular $E$ is slope $H_1$-semistable. Consider the Harder-Narasimhan filtration relative to $H$
$$0 \subset F_1 \subset F_2 \subset \dots \subset F_s=E.$$ By definition, the reduced $H$-Hilbert polynomials of $F_i/F_{i-1}$ are strictly decreasing. In particular, the $H$-slopes are decreasing as well. In turn, this implies $$\mu_{H_1}(F_1)\geq \mu_{H_1}(F_2/F_1)\geq \ldots \geq \mu_{H_1}(F_s/F_{s-1}),$$ and therefore $$\mu_{H_1}(F_1)\geq \mu_{H_1}(F_2)\geq \ldots \geq \mu_{H_1}(F_s)=\mu_{H_1}(E).$$ Since $E$ is slope $H_1$-semistable, we must have
equality throughout $$\mu_{H_1} (F_1) = \mu_{H_1} (F_2) = \ldots = \mu_{H_1} (E).$$ Equivalently, writing $$v(F_i/F_{i-1})=v_i\text{ so that } v=\sum_{i=1}^{s} v_i,$$ we obtain 
\begin{equation}
\label{slopeequal}
\frac{c_1(v_i) \cdot H_1}{\text{rk }v_i}=\frac{c_1(v)\cdot H_1}{\text{rk } v}, \, \,\, \, 1 \leq i \leq s.
\end{equation}

Let ${\mathcal F}_{H} (v_1,v_2,\dots,v_s)$ be the stack of the Harder-Narashimhan
filtrations
\begin{equation}
0 \subset F_1 \subset F_2 \subset \dots \subset F_s=E, \, \, \, E \in \ms(v)
\end{equation}
such that the quotients
 $F_i/F_{i-1}$, $1 \leq i \leq s$ are semistable with respect to $H$ and 
 \begin{equation}
 v(F_i/F_{i-1})=v_i.
 \label{vi}
 \end{equation} Thus
$$
\ms_{H_1}(v)^{\mu\text{-}ss} \setminus \ms_H(v)^{ss}
=\cup_{v_1,...,v_s}
{\mathcal F}_H(v_1,v_2,\dots,v_s),
$$
where \eqref{slopeequal} is satisfied.
Then Lemma $5.3$ in \cite{KY} implies 
\begin{equation}\label{kye}
 \dim {\mathcal F}_{H}(v_1,v_2,\dots,v_s)=\sum_{i=1}^s \dim \ms_H(v_i)^{ss}
 +\sum_{i<j} \langle v_i,v_j \rangle.
\end{equation} 

Write $v_i=\ell_iv_i'$ where $v_i'$ is a primitive Mukai vector. It is shown in Proposition $4.3.4$ of \cite {MMY2} that for all $i, j$ we have $$\langle v_i', v_j'\rangle \geq 3$$ unless either $v_i'$ or $v_j'$ is isotropic, and in this case $\langle v_i', v_j'\rangle \geq 1$. We estimate \begin{eqnarray*}(\langle v, v\rangle+1)&-&\dim (\ms_{H_1}(v)^{\mu\text{-}ss} \setminus \ms_H(v)^{ss})\\ &=& (\langle v, v\rangle+1)-\sum_{i<j} \langle v_i, v_j\rangle-\sum_{i=1}^s \dim \ms_H(v_i)^{ss}\\&=&\sum_{i>j} \langle v_i, v_j\rangle -\sum_{i=1}^{s}\left(\dim \ms_H(v_i)^{ss}-\langle v_i, v_i\rangle\right)+1\\&\geq& \sum_{i>j} \ell_i\ell_j \langle v_i', v_j'\rangle - \sum_{i=1}^{s} \ell_i+1\geq \sum_{i>j} \ell_i \ell_j - \sum_{i} \ell_i+1\geq 2.\end{eqnarray*} Indeed, the above inequality is satisfied for $s\geq 4$. The cases $s=2$ and $s=3$ need to be considered separately. The detailed analysis is contained in Proposition $4.3.4$ of \cite{MMY2}. The only possible exceptions correspond to 
\begin {itemize}
\item [-] $s=2$, $\ell_1=1$, $\ell_2=\ell$, $v_2'$ isotropic, $\langle v_1', v_2'\rangle=1$;
\item [-] $s=2$, $\ell_1=1$, $\ell_2=1$, $v_1'$ isotropic, $\langle v_1', v_2'\rangle=2$;
\item [-] $s=3$, $\ell_1=\ell_2=\ell_3=1$, $v=v_1'+v_2'+v_3'$, $v_i'$ isotropic, $\langle v_i', v_j'\rangle =1$.
\end {itemize}
In all cases, taking $u=v_1'$, we obtain $\langle v, u\rangle =1 \text{ or }2$, which contradicts our assumption. 

For the final claim about the moduli space $\mathfrak K_H(v)^{ss}$, we repeat the proof above. The only modification is the dimension estimate \eqref{kye} which follows by going over the argument in \cite {KY}. 
\qed

\begin {lemma} Assume that $$\langle v, v\rangle> 4\text{ rank } (v).$$ Then no isotropic vector $u$ of positive rank satisfying $\langle v, u\rangle =1\text { or } 2$ occurs as Mukai vector of a quotient in a Harder-Narasimhan filtration of a sheaf of type $v.$ Therefore, the moduli spaces $$\ms_{H}(v)^{ss} \text { and } \mathfrak K_H(v)^{ss}$$ are independent of the polarization $H$ in codimension $1$. 
\end {lemma} 
\proof Assume that there exists an isotropic vector $u$ as above such that $\langle v, u\rangle =1 \text{ or } 2$. In this situation, we have $$\frac{c_1(u)\cdot H_1}{\text{rk }u}=\frac{c_1(v)\cdot H_1}{\text{rk }v}\implies \left(\frac{c_1(u)}{\text{rk }u}-\frac{c_1(v)}{\text{rk} v}\right)\cdot H_1=0.$$
Using the Hodge index theorem, we conclude that $$\left(\frac{c_1(u)}{\text{rk }u}-\frac{c_1(v)}{\text{rk} v}\right)^2\leq 0.$$ By direct calculation, or via Lemma $1.1$ of \cite {KY}, we obtain \begin{eqnarray*}\langle v, u\rangle &=& -\frac{\text{rk} (v)\cdot \text{rk} (u)}{2}\left(\frac{c_1(u)}{\text{rk }u}-\frac{c_1(v)}{\text{rk} v}\right)^2+\frac{\text{rk }(u)}{\text{rk }(v)}\cdot \frac{\langle v, v\rangle}{2}+\frac{\text{rk }{(v)}}{\text{rk }(u)}\cdot \frac{\langle u, u\rangle}{2}\\ & \geq& \frac{\text{rk }(u)}{\text{rk }(v)}\cdot \frac{\langle v, v\rangle}{2}> 2\text{ rk } (u)\geq 2.\end{eqnarray*} This contradiction completes the proof. \qed

\begin {remark} The Lemma above applies to the particular situation of a product abelian surface $X=B\times F$ considered in Section \ref{deg}. We assume here that $B, F$ are not isogenous, so that the section $\sigma$ and the fiber class $f$ generate the N\'eron-Severi group. Then, for Mukai vectors $$v=(r, \sigma+n f, \chi), \,\, w=(r', \sigma+nf, \chi')$$ with $\chi, \chi'<0$ we obtain $$\langle v, v\rangle =2n-2r\chi>-2r\chi\geq 4r,$$ as required in order to apply the Lemma. 

The only exception may be the case $\chi=-1$ which will be treated separately. In this situation, we claim that there are no walls between the polarizations $$H=\sigma+nf, \,\,\,\,H'=\sigma+Nf,$$ where $N$ is taken sufficiently large to ensure that $H'$ is suitable. Indeed, assuming otherwise, consider a wall defined by an isotropic Mukai vector $u$ such that $$\langle v, u\rangle =1 \text { or } 2.$$ In fact, possibly doubling $u$, it suffices to analyze the case $\langle v, u\rangle =2$. Let $$H_0=\sigma+kf,\,\,\,k\geq n,$$ be an ample divisor on this wall, where $k\in \mathbb Q$. By definition, the vector $u$ appears as the Mukai vector of a quotient in the Harder-Narasimhan filtration for $H_0$. Setting $u=(p, \eta, q)$ with $p>0$, we obtain from \eqref{slopeequal} that $$\left(\frac{\eta}{p}-\frac{H}{r}\right)\cdot H_0=0\implies (r\eta-pH)\cdot H_0=0.$$ Writing $$r\eta-pH=a\sigma+bf,$$ we calculate $$(r\eta-pH)\cdot H_0=(a\sigma+bf)\cdot (\sigma+kf)=0\implies b=-ak.$$ Consequently, \begin{equation}\label{sqr}(r\eta-pH)^2=(a\sigma+bf)^2=2ab=-2a^2k\leq -2k,\end{equation} unless $a=b=0$. This particular situation can be analyzed by exactly the same methods; we leave the verification to the reader. In any case, the conditions that $u$ is isotropic and $\langle v, u\rangle =2$ translate into $$\eta^2=2pq,\,\,\, \eta\cdot H=-p+qr+2,$$ respectively. With this understood, we compute the left hand side of \eqref{sqr} $$(r\eta-pH)^2=r^2\eta^2+p^2H^2-2pr(\eta\cdot H)=2np^2+2pr(p-2)\geq 0>-2k$$ with the only possible exception $p=1$. In this case, the above calculation yields $$(r\eta-pH)^2=2n-2r.$$ By orthogonality, $$2n=-r'\chi-r\chi'\geq r+r'>r$$ which implies $$(r\eta-pH)^2=2n-2r>-2n\geq -2k.$$ This contradicts \eqref{sqr}, showing that there is no wall separating $H$ from a suitable polarization. 


\end {remark} \qed

\section{The Verlinde sheaves are locally free}
\label{locallyfree}
The goal of this section is to prove Theorem \ref{t4}. We show that for any $(X, H)$, the dimension of the space of sections of the theta line bundles is given by the expected formula \eqref{h1} for a very general class of Mukai vectors. This holds even without knowing the vanishing of higher cohomology. 
As a consequence, the Verlinde sheaves $\mathbb V$ and $\mathbb W$ used in the degeneration argument of Section \ref{deg} are in fact locally free over the entire moduli space $\mathcal A$ of pairs $(X, H)$. 

The result should be compared to Proposition \ref{bignef} of Section \ref{deg}. The {\it generic} local-freeness yielded by Proposition \ref{bignef} was sufficient for proving our main Theorem \ref{t1}. By contrast, Theorem 4 gives {\it global} local-freeness in great generality, and will be useful for future strange duality studies.  

We split the theorem into two statements with proofs of different flavors. First, we show

\begin {proposition} \label{nohc} Let $(X, H)$ be a polarized abelian surface. Assume that $$v=(r, dH, \chi), \,\,w=(r', d'H, \chi')$$ are orthogonal primitive Mukai vectors of ranks $r, r'\geq 2$ such that \begin {itemize}
\item [(i)] $d, d'>0$;
\item [(ii)] $\chi<0,\,\, \chi'<0.$
\end {itemize} Assume furthermore that if $(d,\chi)=(1, -1)$, then $(X, H)$ is not a product of two elliptic curves. We have $$h^0(\mathsf K_v, \Theta_w)=\chi(\mathsf K_v, \Theta_w)=\frac{d_v^2}{d_v+d_w}\binom{d_v+d_w}{d_v}.$$ 
\end {proposition}

In the same context, the Proposition implies the requisite statement for the moduli space $\mathsf M_v$: 

\begin{proposition}\label{l2}
In the setup of Proposition \ref{nohc}, for any representative $F\in \mathsf K_w$ we have $$h^0(\mathsf M_v, \Theta_F)=\chi(\mathsf M_v, \Theta_F)=\frac{d_w^2}{d_v+d_w}\binom{d_v+d_w}{d_v}.$$ 
\end{proposition} 
\vskip.1in

\subsection{Proof of Proposition \ref{nohc}.} We begin by explaining the strategy of the proof when $\mathsf K_v$ is smooth. The key point is Lemma \ref{l3} below which shows that $\Theta_w\to \mathsf K_v$ is movable, hence (big and) nef on a smooth birational model of $\mathsf K_v$, cf. Theorem 7 of \cite{HT}. The birational models of $\mathsf K_v$ arise as moduli spaces of Bridgeland stable objects. The dimension calculation is carried out on the moduli space of Bridgeland stable objects, where the higher cohomology vanishes. The Proposition follows since wall-crossings do not change the dimension of the space of sections. The case when $\mathsf K_v$ may be singular requires first to desingularize the moduli space. The above argument can then be repeated on a symplectic resolution. 

Let us elaborate the discussion. As already remarked, the proof uses moduli spaces of Bridgeland stable objects. Specifically, we consider stability conditions $\sigma=\sigma_{s, t}=(Z_{s, t}, \mathcal A_{s, t}),$ for $t>0$, corresponding to central charges $$Z_{s, t}(E)=\langle \exp((s+it)H), v(E)\rangle.$$ The heart $\mathcal A_{s, t}$ has as objects certain $2$-step complexes, and is obtained as a tilt of the abelian category of coherent sheaves on $X$ at a certain torsion pair; the exact definition will not be used below, but we refer the reader to \cite {Br} for details. We form the moduli spaces $\mathsf M_v(\sigma)$ of $\sigma$-semistable objects of type $v$. The moduli space comes equipped with the Albanese map $$\mathsf a: \mathsf M_v(\sigma)\to X\times \widehat X,$$ and we write $\mathsf K_v(\sigma)$ for the Albanese fiber. 
\vskip.1in
We begin by analyzing the case $\mathsf K_v$ smooth. The following observations (a)-(c) are useful for the argument. 
\begin {itemize}
\item [(a)] In the large volume limit $t>>0,$ Bridgeland stability with respect to $\sigma_{s, \infty}:=\sigma_{s, t}$ coincides with Gieseker stability, cf. \cite {Br}, Section $14$. 
\end {itemize} 

The next remarks (b)-(c) are contained in the recent papers \cite {MMY} and \cite {Y2}. For $K3$ surfaces, the similar statements are found in \cite {BM}. 
\begin {itemize} 
\item [(b)] The space of stability conditions admits a wall and chamber decomposition, so that the moduli spaces are constant in each chamber, but they undergo explicit birational transformations as walls are crossed. These birational transformations are regular in codimension $1$.  
\end {itemize}

For the next remark, observe that the theta map \eqref{thetamap} gives an isomorphism $$\Theta: (v^{\vee})^{\perp}\to \text{Pic} (\mathsf K_v(\sigma)),$$ in such a fashion that the Beauville-Bogomolov form on the right hand side corresponds to the Mukai pairing on the left hand side.
Two basic (real) cones of divisors are necessary for our purposes. First, the positive cone $$\text{Pos}(\mathsf K_v(\sigma))\hookrightarrow \text{Pic}(\mathsf K_v(\sigma))_{\mathbb R}$$ can be expressed via the Beauville-Bogomolov form $$\text{Pos}(\mathsf K_{v}(\sigma))=\{x: \langle x, x\rangle>0,\,\, \langle x, A\rangle >0 \text{ for a fixed ample divisor } A \text{ over } \mathsf K_v(\sigma)\}.$$
Second, the movable cone $$\text{Mov}(\mathsf K_v(\sigma))\hookrightarrow \text{Pic}(\mathsf K_v(\sigma))_{\mathbb R}$$ is generated by divisors whose stable base locus has codimension $2$ or higher. Positive movable divisors are big and nef on some smooth birational models, cf. Theorem $7$ of \cite {HT}. In our context, we have the following result obtained via the study of the movable cone in \cite {MMY2}:
\begin {itemize}
\item [(c)] A positive movable divisor $$\text{Mov}(\mathsf K_v(\sigma_{s, \infty}))\cap \text{Pos}(\mathsf K_v(\sigma_{s, \infty}))$$ is identified, under the birational wall crossings of (b), with a big and nef divisor on a smooth moduli space $\mathsf K_v(\sigma_{s, t})$ of Bridgeland stable objects:
$$\{ \Theta_w \to   \mathsf K_{v}(\sigma_{s, \infty}) \}\, \, \, \, \, \longleftrightarrow \, \, \, \, \, \{\Theta_{w} \to \mathsf K_{v}(\sigma_{s, t}) \}.$$ (Note that the Mukai vector $w$ labeling the theta line bundle may undergo Weyl reflections when crossing divisorial walls in $(v^{\vee})^{\perp}$. However, since there are no divisorial walls within the movable chamber, $w$ does not change in the present setting.)

\end {itemize} 

The essential ingredient is then provided by the following

\begin{lemma}\label{l3} For $v$ and $w$ as in Proposition \ref{nohc}, the line bundle $\Theta_w\to \mathsf K_v(\sigma_{s, \infty})$ belongs to the positive movable cone.  
\end{lemma}

As a consequence of remarks (a)-(c) and of the lemma, we note
$$h^0(\mathsf K_v(\sigma_{s, \infty}), \Theta_w) = h^0(\mathsf K_v(\sigma_{s, t}), \Theta_{w})=\chi(\mathsf K_{v}(\sigma_{s, t}), \Theta_{w}).$$
By the same argument as for the usual Gieseker stability, as in Proposition $1$ of \cite {abelian1}, we further have
$$\chi(\mathsf K_{v}(\sigma_{s, t}), \Theta_{w})=\frac{d_v^2}{d_v+d_w}\binom{d_v+d_w}{d_v}.$$ 
 We conclude that  $$h^0(\mathsf K_v(\sigma_{s, \infty}), \Theta_w)=\frac{d_v^2}{d_v+d_w}\binom{d_v+d_w}{d_v},$$ as claimed in Proposition \ref{nohc}. \qed

\vskip.1in

{\it Proof of Lemma \ref{l3}.} We begin by noting that $\Theta_w$ is positive in the Gieseker chamber. Indeed, $$\langle \Theta_w, \Theta_w\rangle=\langle w, w\rangle> 0.$$ For the second inequality, an ample divisor on the moduli space $\mathsf K_v(\sigma_{s, \infty})$ is constructed in \cite {lepotier}; see also Remark 8.1.12 of \cite {HL}. This divisor takes the form $\Theta_a$ for $$a=(r, rm H, -2mnd-\chi), \, \, \, \, \text{where} \,\, \, m>>0.$$ Recalling that $w=(r', d'H, \chi')$, we have $$\langle \Theta_w, \Theta_a\rangle=\langle w, a\rangle  = 2nm(d'r+dr')-r\chi'+r'\chi>0,$$ as needed. 

We will now show that for the vector $$w_1=(2nd, -\chi H, 0)$$ the line bundle $\Theta_{w_1}$ belongs to the closure of the movable cone for the Gieseker chamber. We will combine this with a well-known result of Jun Li \cite {junli2}. For the vector $$w_0=(0, rH, -2nd)$$ the associated theta line bundle $$\Theta_{w_0}\to \mathsf K_v(\sigma_{s, \infty})$$ is big and nef, so in particular it is in the closure of the movable cone. Notice now that the vector $w$ is a positive linear combination of $w_0$ and $w_1$, 
$$w = \frac{1}{2nd} \left ( -\chi' w_0 + r' w_1 \right ),$$
hence $\Theta_w$ is movable. 

To prove the claim about $w_1$, we will use the description of the movable cone given in \cite {MMY} and \cite {BM}. Specifically, we consider the hyperplanes in $\overline{\text{Pos}}(\mathsf K_{v}(\sigma_{s, \infty}))$ given by $$\Theta((u^\vee)^{\perp}\cap (v^{\vee})^{\perp}),\,\,\, 1\leq \langle v, u\rangle \leq 2, \,\, \langle u, u\rangle=0.$$ The movable cone is cut out by these hyperplanes.  To prove that $\Theta_{w_1}$ and $\Theta_{w_0}$ belong to the same chamber, it suffices to show that $$\langle w_0, u^{\vee}\rangle \geq 0 \iff \langle w_1, u^{\vee}\rangle \geq 0,$$ whenever $u$ is isotropic and $1\leq \langle v, u\rangle \leq 2.$ The first inequality above will in fact turn out strict for rank $3$ or higher. 

We assume $r>2$ first. Let us write $u=(p, \eta, q)$ where $${\eta^2}= 2pq, \, \, \, \, p, q \in {\mathbb Z}.$$ Changing $u$ into $-u$, we may furthermore assume that $p\geq 0$ and $\langle v, u\rangle =\pm 1, \pm 2$. Recalling that $v=(r, dH, \chi)$, we calculate \begin{equation}\label{par}\langle v, u\rangle = d(H\cdot \eta)-p\chi-qr=\pm 1\text { or } \pm 2.\end{equation} We compute \begin{equation}\label{e1}\langle w_0, u^{\vee}\rangle \geq 0\iff -r(H\cdot \eta)+2ndp\geq 0.\end{equation} Similarly, \begin{equation}\label{e2}\langle w_1, u^{\vee}\rangle \geq 0\iff \chi (H\cdot \eta)-2ndq\geq 0.\end{equation} We therefore need to show that 
$$-r(H\cdot \eta)+2ndp\geq 0 \iff \chi (H\cdot \eta)-2ndq\geq 0.$$

We consider first the case when $p=0$. Then, replacing $u$ by $-u$ we may assume that $H\cdot \eta \geq 0$. In fact, $H\cdot \eta=0$ is impossible by \eqref{par} since $r>2$. Therefore, $H\cdot \eta>0$. In this situation, \eqref{e1} is false. We argue that \eqref{e2} is false as well. Assuming otherwise, we have $$\chi (H\cdot \eta)\geq 2ndq\implies q< 0.$$ This is however incompatible with \eqref{par} which reads $$d(H\cdot \eta)+r(-q)=\pm 1, \pm 2,$$ which is impossible for $r>2$.

The crux of the argument is the case $p>0$. In this situation, we distinguish the following subcases:
\begin {itemize}
\item [(i)] Assume $H\cdot \eta=0$. By the Hodge index theorem $\eta^2\leq 0$ hence $$pq=\frac{\eta^2}{2}\leq 0\implies q\leq 0.$$ This shows that both \eqref{e1} and \eqref{e2} are true at the same time. 
\item [(ii)] Assume $H\cdot \eta<0$. In this case, \eqref{e1} is true. We prove that \eqref{e2} is true as well. Assuming otherwise, we obtain that $$\chi (H\cdot \eta)-2ndq< 0.$$ In particular $q>0$ and multiplying by $p>0$ we see that $$\frac{p\chi}{2nd} (H\cdot \eta)<pq=\frac{\eta^2}{2}.$$ By the Hodge index theorem, we have $$\eta^2\leq \frac{(H\cdot \eta)^2}{2n}.$$ The above inequality becomes $$\frac{p\chi}{2nd} (H\cdot \eta)<\frac{(H\cdot \eta)^2}{4n}\implies (H\cdot \eta)< \frac{2p\chi}{d}.$$ We obtain therefore 
$$d(H\cdot \eta)-p\chi-qr< 2p\chi-p\chi-qr=p\chi-qr < -2,$$ using $\chi<0$ and $q>0$. This contradicts \eqref{par}. Thus \eqref{e2} must be true as well. 
\item [(iii)] Assume $H\cdot \eta>0$. Equation \eqref{par} implies that $q\geq 0$. In this case, the inequality \eqref{e2} is false. We argue that \eqref{e1} is false as well. Assume otherwise, so that $$r(H\cdot \eta)\leq 2ndp\implies \frac{rq}{2nd}(H\cdot \eta)\leq pq=\frac{\eta^2}{2}.$$ Again by the Hodge index theorem, we have $$\eta^2\leq \frac{(H\cdot \eta)^2}{2n}$$ yielding $$\frac{rq}{2nd}(H\cdot \eta)\leq pq=\frac{\eta^2}{2}\leq \frac{(H\cdot \eta)^2}{4n}\implies 2rq\leq d(H\cdot \eta).$$ We obtain $$d(H\cdot \eta)-p\chi-rq\geq 2rq-p\chi-rq=rq-p\chi> 2$$ if $q>0$, contradicting \eqref{par}. When $q=0,$ equation \eqref{par} yields $$d(H\cdot \eta)-p\chi=\pm 1, \pm 2,$$ which implies $d(H\cdot \eta)=1, p\chi=-1$. Therefore $(d, \chi)=(1, -1)$ and $H\cdot \eta=1$, $\eta^2=0$. In this case, $(X, H)$ is a product of elliptic curves, which is not allowed.\footnote{To see that $X$ is a product, write $\tau=H-n\cdot \eta.$ Therefore, $$\eta^2=\tau^2=0,\,\,\, \eta\cdot \tau=1.$$ In this situation, $\eta$ and $\tau$ are represented by two elliptic curves $E$ and $F$, cf. Proposition $2.3$ in \cite {K}. The sum morphism $$s: E\times F\to X$$ must be an isogeny. The preimage of the origin corresponds to the intersection $E\cap F$, hence $s$ must be an isomorphism. 
} \end {itemize} 

When $r=2$, the same argument goes through with the only exception corresponding to the case $$p=0, \,\,H\cdot \eta=0.$$ Since $\eta^2=2pq=0$ we obtain $\eta=0$ by the Hodge index theorem. This yields the isotropic vector $u=(0, 0, 1)$. In fact, $w_0$ lies on the wall determined by $u$, hence we cannot pin down on which side of the wall $w_1$ lies. To remedy this problem, we replace $w_0$ by the vector $$a=(r, rmH, -2mnd-\chi)=mw_0+(r, 0, -\chi)$$ which we have already seen to give an ample theta bundle for $m>>0$. For the vector $u=(0, 0, 1)$, direct computation shows $$\langle a, u\rangle < 0,\,\, \langle w_1, u\rangle <0,$$ hence $w_1$ and $a$ are also on the same side of the wall determined by $u$. 

This completes the analysis, and therefore the proof when $\mathsf K_v$ is smooth. \vskip.1in

However, $\mathsf K_v$ may be singular when the polarization $H$ is not generic. In this situation, for any $\beta\in \text{NS}(X)_{\mathbb Q}$, we consider the moduli space of $\beta$-twisted $H$-semistable sheaves. Recall that a sheaf $E$ is $\beta$-twisted $H$-semistable provided that \begin{itemize}
\item [(i)] for all subsheaves $F\subset E$, we have $$\frac{c_1(F)\cdot H}{\text{rk}(F)}\leq \frac{c_1(E)\cdot H}{\text{rk}(E)};$$
\item [(ii)] if equality holds in (i), then $$\frac{\chi(F)-c_1(F)\cdot \beta}{\text{rk}(F)}\leq\frac{\chi(E)-c_1(E)\cdot \beta}{\text{rk}(E)}.$$
\end {itemize} We form the moduli space $\mathsf K_{\beta}(v)$ of $\beta$-twisted $H$-semistable sheaves. In fact, remark (a) above applies here as well, and consequently, $\mathsf K_\beta(v)$ can be viewed as a moduli space of Bridgeland's stable objects. In addition, if $\beta$ is appropriately chosen, then 
the moduli space $\mathsf K_{\beta}(v)$ consists of stable sheaves only, and therefore is a smooth non-empty holomorphic symplectic manifold; see for instance Lemma $5.4$ of \cite {abe}. Furthermore, Lemma $5.5$ in \cite {abe} shows that there is a surjective morphism $$\pi:\mathsf K_{\beta}(v)\to \mathsf K_v,$$ which is therefore a symplectic resolution. As a consequence of Proposition $1.3$ of \cite {beauville2} we have $${\mathbf R}\pi_{\star} \mathcal O_{\mathsf K_{\beta}(v)}=\mathcal O_{\mathsf K_v}.$$ Now, as the moduli space $\mathsf K_\beta(v)$ consists of stable sheaves only, it carries a theta line bundle $\Theta_w$. Furthermore, the line bundle $\Theta_w$ descends to the singular moduli space $\mathsf K_v$, which may contain strictly semistables. This is a consequence of Kempf's lemma and is shown to hold true in Theorem $8.1.5$ of \cite {HL}. The essential point is that $c_1(v)=dH$ is a multiple of the polarization. As a corollary, $$H^0(\mathsf K_v, \Theta_w)=H^0(\mathsf K_\beta(v), \pi^{\star}\Theta_w)=H^0(\mathsf K_\beta(v), \Theta_w).$$ We claim that $\Theta_w$ is movable over the smooth moduli space $\mathsf K_{\beta}(v)$. In fact, the argument we presented in the untwisted case carries over to the twisted situation. An essential ingredient of the proof is that Jun Li's line bundle is big and nef. This continues to hold over $\mathsf K_{\beta}(v)$ by pullback, at least for $\beta$ chosen as above. Alternatively, ample divisors are constructed in Lemma $5.5.2$ of \cite {MMY}. Since $\Theta_w$ is movable, we conclude that $$h^0(\mathsf K_\beta(v), \Theta_w)=\chi(\mathsf K_\beta(v), \Theta_w)=\frac{d_v^2}{d_v+d_w}\binom{d_v+d_w}{d_v},$$ as claimed. This completes the proof. 

\qed
\begin {remark}
The argument above also remains valid in ranks $0$ and $1$. Consequently, the dimension calculation of Proposition \ref{nohc} holds true for all primitive orthogonal Mukai vectors $$v=(r, dH, \chi), \,\,w=(r', d'H, \chi') \text{ with }r, r'\geq 0, \,\,\, d, d'>0, \,\,\, \chi, \chi'<0,$$ with the extra assumption that 
\begin {itemize} 
\item [-] $(X, H)$ is not a product when $(d, \chi)=(1, -1)$ or when $(r, d)=(1, 1)$. 
\end {itemize} 
\end {remark} 

\subsection {Proof of Proposition \ref{l2}.} To prove the Proposition, we use the diagram \begin{center} $\xymatrix{\mathsf K_{v}\times X \times \widehat
X\ar[r]^{\hskip.2in\Phi_v} \ar[d]^{p} & \mathsf M_v \ar[d]^{\mathsf a} \\ X\times \widehat X \ar[r]^{\Psi_v} & X\times \widehat X}.$ \end{center} Here,  
$\Phi_v:\mathsf K_{v}\times
X\times \widehat X\to \mathsf M_v$ is defined as $$\Phi_v(E, x, y)=t_{x}^{\star}E \otimes y,$$ and $$\mathsf a: \mathsf M_v\to X\times \widehat X$$ is the Albanese map.  Both $\Phi_v$ and $\Psi_v$ are \'etale of degree $d_v^4$ \cite {Y}, \cite {abelian1}. In fact, it is proved in \cite{Y} that $$\Psi_v(x, y)=(-\chi x- d\varphi_{\widehat H}(y), \,d\varphi_{H}(x)+ry),$$ where as usual $$\widehat H\to \widehat X$$ is the inverse determinant of the Fourier-Mukai transform of $H$, and $\varphi_{H}$, $\varphi_{\widehat H}$ denote the Mumford homomorphisms. This explicit expression will however not be needed below.

Fix $F\in \mathsf K_w$. We have $$\Phi_v^{\star} \Theta_F=\Theta_w\boxtimes \mathcal L$$ for a line bundle $\mathcal L\to X\times \widehat X.$ It is shown in Proposition $4$ of \cite {abelian1} that $$\chi(\mathcal L)=d_v^2d_w^2.$$ In fact, by Lemma $1$ in \cite {bundle}, up to numerical equivalence we have 
\begin{equation}
\label{abc}
\mathcal L=H^{a}\boxtimes \widehat {H}^{b}\otimes \mathcal P^c,
\end{equation}
where $\mathcal P\to X\times \widehat X$ is the Poincar\'e bundle, and $$a=-(\chi d'+\chi'd),\,\,\, b=rd'+r'd,\,\,c=dd'n+r'\chi=-dd'n-r\chi'.$$ In consequence of the assumptions $\chi, \chi'<0$ and $d, d'>0$, and also of the calculation $$abn-c^2=d_vd_w>0,$$  we obtain the inequalities
$$a>0,\,\,\, b>0,\,\,\,abn>c^2.$$ 
These inequalities ensure that the line bundle $\mathcal L$ is ample. To see this, we use the special form of the Nakai-Moishezon criterion for ampleness in the context of abelian varieties, as stated  on page 77 of \cite{bl}. Specifically, for abelian varieties, the criterion asserts that it is enough to check ampleness numerically on hyperplanes and intersections of hyperplanes under any fixed projective embedding, such as the one induced by $H+\widehat H$. A direct calculation then shows that a line bundle $\mathcal L\to X \times \widehat X$ of the form \eqref{abc} is ample if and only if the three inequalities above are satisfied. In consequence, $\mathcal L$ has no higher cohomology. 

With this understood, we write with the aid of Proposition \ref{nohc} \begin{equation}\label{g1}h^0(\mathsf K_v\times X\times \widehat X, \Theta_w\boxtimes \mathcal L)=h^0(\mathsf K_v, \Theta_w)h^0(X\times \widehat X, \mathcal L)=\chi(\mathsf K_v, \Theta_w)\chi(X\times \widehat X, \mathcal L)\end{equation} $$=\frac{d_v^2}{d_v+d_w}\binom{d_v+d_w}{d_v}\cdot (d_vd_w)^2.$$ On the other hand, \begin{equation}\label{g2}h^0(\mathsf K_v\times X\times \widehat X, \Phi_v^{\star}\Theta_F)=h^0(\mathsf M_v, (\Phi_v)_{\star}\Phi_v^{\star} \Theta_F)=\sum_{\tau} h^0(\mathsf M_v, \Theta_F\otimes \mathsf a^{\star}\mathbb L_\tau)\end{equation} where $$(\Psi_v)_{\star} \mathcal O=\bigoplus_{\tau} \mathbb L_{\tau},$$ over $X\times \widehat X$. The line bundles $\mathbb L_{\tau}$ appearing in the decomposition above are indexed by the characters $\tau\in \widehat {\mathsf G}_v$ of the group $$\mathsf G_v=\text{Ker } \Psi_v.$$
We claim that \begin{lemma}\label{aut} For each character $\tau$ of $\mathsf G_v$, there exists an automorphism $f_{\tau}:\mathsf M_v\to \mathsf M_v$ such that $$\Theta_F\otimes \mathsf a^{\star}\mathbb L_{\tau}=f_{\tau}^{\star}\Theta_F.$$
\end {lemma}  

\noindent By the lemma, we therefore have $$h^0(\mathsf M_v, \Theta_F\otimes \mathsf a^{\star}\mathbb L_{\tau})=h^0(\mathsf M_v, \Theta_F)$$ hence by \eqref{g2} we obtain $$h^0(\mathsf K_v\times X\times \widehat X, \Phi_v^{\star}\Theta_F)=\deg \Psi_v \cdot h^0(\mathsf M_v, \Theta_F)=d_v^4 \cdot h^0(\mathsf M_v, \Theta_F).$$ This implies via \eqref{g1} that $$h^0(\mathsf M_v, \Theta_F)=\frac{d_w^2}{d_v+d_w}\binom{d_v+d_w}{d_v},$$ establishing Proposition \ref{l2}.
\qed

{\it Proof of Lemma \ref{aut}.} We consider the group $$K(\mathcal L)\hookrightarrow X\times \widehat X$$ of pairs $(x, y)$ leaving $\mathcal L$ invariant by translation $$t_{(x,y)}^{\star}\mathcal L\simeq\mathcal L.$$ The group $K(\mathcal L)$ has $\chi(\mathcal L)^2=(d_vd_w)^4$ elements. 

For each pair $(x, y)\in K(\mathcal L)$, we define the automorphism $$f_{(x, y)}:\mathsf M_v\to \mathsf M_v$$ given by $$f_{(x, y)}(E)=t_x^{\star}E\otimes y.$$ We show that for $(x, y)\in K(\mathcal L)$ we can find a line bundle $\mathbb L_{\tau}\in \widehat {\mathsf G}_v$ such that \begin{equation}\label{aut1}f_{(x, y)}^{\star} \Theta_F=\Theta_F\otimes {\mathsf a}^{\star} \mathbb L_{\tau}.\end{equation}
Indeed, the two lines bundles $f_{(x, y)}^{\star} \Theta_F$ and $\Theta_F$ both restrict to $\Theta_w$ on each fiber of the Albanese map ${\mathsf a}$, hence for some line bundle $\mathbb L\to X\times \widehat X$ we have $$f_{(x, y)}^{\star} \Theta_F=\Theta_F\otimes {\mathsf a}^{\star} \mathbb L.$$ It remains to explain that $$\Psi_v^{\star}\,\mathbb L=\mathcal O,$$ or equivalently that $$\Phi_v^{\star} f_{(x, y)}^{\star} \Theta_F=\Phi_v^{\star} \Theta_F.$$ Direct calculation shows that over $\mathsf K_v\times X\times \widehat X$ we have $$f_{(x, y)}\circ \Phi_v=\Phi_v\circ (1, t_{(x, y)}).$$ 
Therefore \begin{eqnarray*}\Phi_v^{\star} f_{(x, y)}^{\star}\Theta_F&=&(1, t_{(x,y)})^{\star} \Phi_v^{\star} \Theta_F=(1, t_{(x,y)})^{\star} (\Theta_w\boxtimes \mathcal L)\\ &=&\Theta_w\boxtimes t_{(x, y)}^{\star}\mathcal L=\Theta_w\boxtimes \mathcal L=\Phi^{\star}_v \Theta_F.\end{eqnarray*}

As a consequence of \eqref{aut1}, there exists a group homomorphism $$\alpha:K(\mathcal L)\to \widehat {\mathsf G}_v.$$ To complete the proof of the Lemma, we argue that $\alpha$ is surjective. 
Since $$\text{order } K(\mathcal L)=(d_vd_w)^4, \, \text {order } \mathsf G_v=d_v^4$$ it suffices to prove that $$\text{order } \text{Ker }\alpha=d_w^4.$$ In fact, we claim that \begin{equation}\label{kernel}\text{Ker } \alpha \simeq \mathsf G_w,\end{equation} where $\mathsf G_w$ is the kernel of the morphism $\Psi_w$ in the diagram
\begin{center} $\xymatrix{\mathsf K_{w}\times X \times \widehat
X\ar[r]^{\hskip.2in\Phi_w} \ar[d]^{p} & \mathsf M_w \ar[d]^{\mathsf a} \\ X\times \widehat X \ar[r]^{\Psi_w} & X\times \widehat X}.$ \end{center} Here,  
$\Phi_w: \mathsf K_{w}\times
X\times \widehat X\to \mathsf M_w$ is defined as $$\Phi_w(G, x, y)=t_{-x}^{\star}G \otimes y,$$ and $$\mathsf a: \mathsf M_w\to X\times \widehat X$$ is the Albanese map $$\mathsf a(G)=(\det \widehat G\otimes \widehat H^{d'}, \det G\otimes H^{-d'}).$$  
Furthermore, just as above, $\Phi_w$ and $\Psi_w$ both have degree $d_w^4$. To prove \eqref{kernel}, note that $$(x,\, y)\in \text{Ker } \alpha \iff f_{(x, \,y)}^{\star} \Theta_F=\Theta_F\iff \Theta_{t_{-x}^{\star} F\otimes y}=\Theta_F.$$ By \cite {abelian1}, the last equality happens if and only if $$\det (t_{-x}^{\star} F\otimes y)=\det F \text { and } \det \widehat{(t_{-x}^{\star} F\otimes y)}=\det \widehat F$$ $$\iff (\mathsf a\circ \Phi_w)(F, x, y)=0\iff \Psi_w(x, y)=0\iff (x, y)\in \mathsf G_w,$$ as claimed. The proof of the lemma is completed. \qed

\end {document}